\documentclass[oneside,english]{amsart}
 \usepackage[T1]{fontenc}
\usepackage{a4}
\usepackage{graphicx}
\usepackage{setspace}
\usepackage{amssymb}

\makeatletter
 \theoremstyle{plain}    
 \newtheorem{thm}{Theorem}[section]
 \numberwithin{equation}{section} %% Comment out for sequentially-numbered
 \numberwithin{figure}{section} %% Comment out for sequentially-numbered
 \theoremstyle{plain}
 \theoremstyle{remark}
 \newtheorem*{rem*}{Remark}
 \theoremstyle{definition}
  \newtheorem*{example*}{Example}
 \theoremstyle{definition}
 \newtheorem{defn}[thm]{Definition}
 \theoremstyle{plain}    
 \newtheorem{prop}[thm]{Proposition} %%Delete [thm] to re-start numbering
 \theoremstyle{plain}    
 \newtheorem{lem}[thm]{Lemma} %%Delete [thm] to re-start numbering
 \theoremstyle{remark}
 \newtheorem{rem}[thm]{Remark}

\usepackage{setspace}\usepackage{psfrag}
 \usepackage{amssymb}
\DeclareMathOperator*{\diam}{diam}

\DeclareMathOperator*{\Fin}{Fin}
\DeclareMathOperator*{\card}{card}

\DeclareMathOperator*{\HD}{HD}

\DeclareMathOperator*{\Con}{Con}
\DeclareMathOperator*{\Int}{Int}
\DeclareMathOperator*{\sgn}{sign}
\DeclareMathOperator*{\dom}{dom}

\newcommand{\N}{\mathbb{N}}
\renewcommand{\P}{\mathcal{P}}

\renewcommand{\hat}{\widehat}
\newcommand{\K}{\mathcal{K}}

\renewcommand{\epsilon}{\varepsilon}
\newcommand{\R}{\mathbb{R}}

\usepackage{xr-hyper}
\newcommand{\e}{\mathrm{e}}
\usepackage[colorlinks=false,bookmarks,pdfhighlight=/O,pdfauthor={Marc Kesseboehmer},pdftitle={Multidimensional}]{hyperref}

\usepackage{babel}

\usepackage{babel}
\makeatother
\begin{document}
\title[Higher-dimensional multifractal value sets for infinite systems]{Higher-dimensional multifractal value sets for conformal infinite graph directed Markov systems}

\author{$^{1}$Marc Kesseböhmer and $^{2}$Mariusz Urba\'{n}ski}

\begin{abstract}
We give a description of the level sets in the higher dimensional
multifractal formalism for infinite conformal graph directed Markov
systems. If these systems possess a certain degree of regularity this
description is complete in the sense that we identify all values with
non-empty level sets and determine their Hausdorff dimension. This
result is also partially new for the finite alphabet case. 
\end{abstract}

\date{\today}

\address{$^{1}$Universität Bremen, Fachbereich 3 - Mathematik und Informatik,
Bibliothekstrasse 1, 28334 Bremen, Germany. Tel.: +49-(0)421-218-8619,
Fax: +49-(0)421-218-4235.}

\email{mhk@math.uni-bremen.de}

\address{$^{2}$University of North Texas, Department of Mathematics, P.O.
Box 311430, Denton, TX 76203, USA.}

\email{urbanski@unt.edu}

\keywords{Multifractal fromalism, thermodynamic formalism, conformal graph
directed Markov systems.}

\subjclass{37C45, 28A80, 82B05}

\maketitle

\section{Introduction and statements of main results}

\let\languagename\relax

In this paper we study \emph{Graph Directed Markov System} (GDMS)
as defined in \cite{MauldinUrbanski:03} consisting of a directed
multigraph $(V,E,i,t,A)$ with incidence matrix $A$ together with
a family of non-empty compact metric spaces $\left(X_{v}\right)_{v\in V}$,
a number $s\in\left(0,1\right)$, and for every $e\in E$, an injective
contraction $\phi_{e}:X_{t(e)}\rightarrow X_{i(e)}$ with Lipschitz
constant not exceeding $s$. Briefly, the family\[
\Phi=\left(\phi_{e}:X_{t(e)}\rightarrow X_{i(e)}\right)_{e\in E}\]
 is called a GDMS. Throughout this paper we will assume that the system
is \emph{conformal} (Def. \ref{def:Conformal}), \emph{finitely irreducible}
(Def. \ref{def:finiteIrreducible}), and \emph{co-finitely regular}
(Def. \ref{d4.2.2}). The necessary details will be postponed to Section
\ref{sec:Conformal-Graph-Directed}. Let $E_{A}^{\infty}$ denote
the set of admissible infinite sequences for $A$ and $\sigma:E_{A}^{\infty}\to E_{A}^{\infty}$
the shift map given by $\left(\sigma\left(x\right)\right)_{i}:=\left(x_{i+1}\right)_{i}$.
With $\pi:E_{A}^{\infty}\to X:=\bigoplus_{v\in V}X_{v}$ we denote
the natural coding map from the subshift $E_{A}^{\infty}$ to the
disjoint union $X$ of the compact sets $X_{v}$ (see (\ref{eq:pi})).
Its image $\Lambda:=\Lambda_{\Phi}:=\pi\left(E_{A}^{\infty}\right)$
denotes the limit set of $\Phi$. An important tool for studying $\Phi$
is the following geometric potential function given by the conformal
derivatives of the contractions $\left(\phi_{i}\right)_{i\in E}$\[
I=I_{\Phi}:E_{A}^{\infty}\rightarrow\R^{+},\: I_{\Phi}\left(\omega\right):=-\log\left|\phi'_{\omega_{1}}\left(\pi\left(\sigma\left(\omega\right)\right)\right)\right|.\]
We are going to set up a multifractal analysis for $I$ with respect
to another $d$-dimensional bounded Hölder continuous function\[
J:E_{A}^{\infty}\to\R^{d}.\]
That is for $v\in V$ and $\alpha\in\R^{d}$ we investigate the level
sets\[
\mathcal{F}_{\alpha}\left(v\right):=\left\{ \pi\left(\omega\right):\omega\in E_{A}^{\infty},i\left(\omega_{1}\right)=v,\,\textrm{and }\lim_{n\rightarrow\infty}\frac{S_{n}J\left(\omega\right)}{S_{n}I\left(\omega\right)}=\alpha\right\} ,\]
and $\mathcal{F}_{\alpha}=\bigoplus_{v\in V}\mathcal{F}_{\alpha}\left(v\right)\subset\Lambda$.
In here, $S_{n}f\left(\omega\right):=\sum_{k=0}^{n-1}f\left(\sigma^{k}\left(\omega\right)\right)$.
Let \[
Q:\mathcal{M}\left(E_{A}^{\infty},\sigma\right)\rightarrow\R^{d},\quad Q\left(\mu\right):=\frac{1}{\int I\, d\mu}\int J\, d\mu,\]
where $\mathcal{M}\left(E_{A}^{\infty},\sigma\right)$ denotes the
set of shift invariant Borel probability measures on $E_{A}^{\infty}$.
We are now interested in the following three subsets of $\R^{d}$.
\begin{eqnarray}
K & := & \left\{ \alpha\in\R^{d}:\mathcal{F}_{\alpha}\neq\emptyset\right\} ,\nonumber \\
L & := & \left\{ Q\left(\mu\right)=\left(\int I\, d\mu\right)^{-1}\int J\, d\mu:\mu\in\mathcal{M}\left(E_{A}^{\infty},\sigma\right)\right\} ,\label{eq:DefKLM}\\
M & := & \nabla\beta\left(\R^{d}\right),\nonumber \end{eqnarray}
 where $\beta:\R^{d}\to\R$ is the convex differentiable function
defined in terms of some pressure function within Proposition \ref{pro:PressureBetaGibbs}.
Since $I>-\log s$, the sets $K,L,M\subset\R^{d}$ are all bounded.
Since $E_{A}^{\infty}$ is finitely irreducible (see Def. \ref{def:finiteIrreducible})
we have for all $v\in V$ \[
K=\left\{ \alpha:\mathcal{F}_{\alpha}\left(v\right)\not=\emptyset\right\} .\]
Our first main theorem relates the three sets in (\ref{eq:DefKLM}).

\begin{thm}
\label{thm:Main1} The set $K$ is compact and we have $\Int L\subset M\subset L$
and $\overline{M}\subset K\subset\overline{L}.$ 
\end{thm}
\begin{rem*}
For the finite alphabet case (i.e. $E$ is a finite set) the inclusion
$K\subset L$ is well-known for the one dimensional situation (i.e.
$d=1$) and equality of $K$ and $L$ is also proved for $d\geq1$
in \cite{BarreiraSaussolSchmeling:02}. The proof uses the fact that
for $x\in\mathcal{F}_{\alpha}\ne\emptyset$ the set of measures $\left\{ \mu_{n}:=n^{-1}\sum_{i=0}^{n-1}\delta_{\sigma^{i}x}:n\in\N\right\} $
always possesses a weak convergent subsequence with limit measure
$\mu$ such that $Q\left(\mu\right)=\alpha$. This gives $\alpha\in L$.
Also, since in the finite alphabet case $\mathcal{M}\left(E_{A}^{\infty},\sigma\right)$
is compact and $Q$ is continuous with respect to the weak{*}-topology
we have $\overline{\Int L}=L$. Hence, the above theorem gives in
this situtation $L=\overline{\Int L}\subset\overline{M}\subset K$.
This shows that $L=K$ in the finite alphabet case. 
\end{rem*}
If some additional regularity conditions are satisfied we get the
following stronger results.

\begin{thm}
\label{thm:Main2} Suppose that $J_{i}$ are linearly independent
as cohomology classes. Then $M$ is an open convex domain, $L\subset\overline{\Int L}$,
and in particular \[
\overline{L}=\overline{M}=K.\]
If additionally $0\in M$ then $L=\overline{M}=K$.
\end{thm}
Our third theorem gives the multifractal formalism in the higher-dimensional
situation. If $0\in M$ then our description is complete in the sense
that the formula for the Hausdorff dimension holds not only for $\alpha$
form the interior of $K$ but for all $\alpha\in\R^{d}$. This result
is also partially new in the finite alphabet case. However, similar
results for the finite alphabet case in the context of mixed singularity
spectra for self-similar measures or deformed Birkhoff averages are
obtained in \cite{Olsen:05Mixed,Olsen:03Mixed,OlsenWinter:07}. Recall
that the (negative) Legendre transform $\hat{\beta}:\R^{d}\to\overline{\R}$
of $\beta$ is given by \[
\hat{\beta}\left(a\right):=\inf_{t\in\R^{d}}\left(\beta\left(t\right)-\left\langle t,a\right\rangle \right).\]
For the following let $\HD\left(A\right)$ denotes the Hausdorff dimension
of the set $A$. 

\begin{thm}
\label{thm:MF} Suppose that $J_{i}$ are linearly independent as
cohomology classes. Then we have for $\alpha\in M$ and $v\in V$\[
\HD\left(\mathcal{F}_{\alpha}\left(v\right)\right)=\widehat{\beta}\left(\alpha\right)\]
and for all $\alpha\in\R^{d}$ we have \begin{equation}
\HD\left(\mathcal{F}_{\alpha}\left(v\right)\right)\leq\max\left\{ \widehat{\beta}\left(\alpha\right),0\right\} .\label{eq:HDupperBound}\end{equation}
If additionally $0\in M$ then equality holds in (\ref{eq:HDupperBound}). 
\end{thm}
\begin{example*}
\begin{figure}[h]
\psfrag{t2}{$t_2$}

\psfrag{beta}{$\beta$}\psfrag{t1t2beta}{ $(t_1,t_2,\beta(t_1,t_2))$}  \psfrag{betahat}{$\widehat{\beta}(\partial_1 \beta(t_1,t_2),\partial_2 \beta(t_1,t_2))$}

\psfrag{t1}{$t_1$}

\includegraphics[%
  width=0.90\columnwidth,
  keepaspectratio]{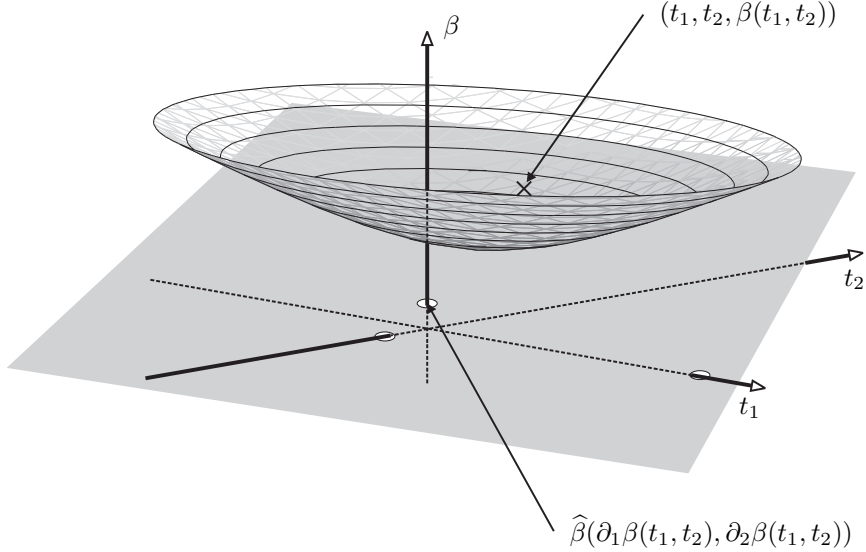}

\caption{The Legendre transform $\hat{\beta}$ at $\nabla\beta(t)$ determined
by the distance between the origin and the unique point at which the
$\beta$-axis intersects the hyperplane tangential to the graph of
$\beta$ at the point $\left(t_{1},t_{2},\beta\left(t_{1},t_{2}\right)\right)$.
\label{cap:fig1}}
\end{figure}
 For an example illustrating the above theorems let us consider the
iterated function system generated by continued fractions, i.e. \begin{equation}
\Phi:=\left(\phi_{k}:\left[1,0\right]\to\left[1,0\right],x\mapsto1/(x+k):k\in\N\right).\label{eq:CF}\end{equation}
 Then $E_{A}^{\infty}=\N^{\N}$ and for $d=2$ let us choose \[
J:\N^{\N}\to\R^{2},\;\left(n_{1},n_{2},...\right)\mapsto J\left(n_{1}\right)\]
with \[
J_{1}\left(n\right):=\left\{ \begin{array}{rll}
-1 & \:\textrm{for} & n=0\mod2\\
1 & \:\textrm{for} & n=1\mod2\end{array},\right.\; J_{2}\left(n\right):=\left\{ \begin{array}{rll}
0 & \:\textrm{for} & n=0\mod3\\
1 & \:\textrm{for} & n=1\mod3\\
-1 & \:\textrm{for} & n=2\mod3\end{array},\right.\]
for $n\in\N$. Now we make use of the fact that $J_{1}$ and $J_{2}$
are linearly independent as cohomology classes if and only if the
following implication holds.\[
\exists C>0\forall\omega\in\N^{\N},\forall n\in\N:\,\left|\alpha_{1}S_{n}J_{1}(\omega)+\alpha_{2}S_{n}J_{2}(\omega)\right|<C\:\implies\:\alpha_{1}=\alpha_{2}=0.\]
 If we consider the periodic infinite words $\omega=\left(3,3,3,\ldots\right)$
and $\left(1,1,1,\ldots\right)$ the above implication follows immediately
for this particular choice of $J$. Moreover, in this situation we
easily verify that for the function $\beta$ as defined in Proposition
\ref{pro:PressureBetaGibbs} we have\[
\lim_{\left|t\right|\to\infty}\beta\left(t\right)=\infty.\]
Therefore $\beta$ has a unique minimum in $\R^{2}$ and consequently
$0\in M$. Hence, we can apply all of the above theorems to deduce
that $L=\overline{M}=K$ and for all $\alpha\in\R^{2}$ we have $\HD\left(\mathcal{F}_{\alpha}\left(v\right)\right)=\widehat{\beta}\left(\alpha\right)$.
Qualitatively, this situation has been illustrated in Figure \ref{cap:fig1}.
\end{example*}
A detailed study of the higher dimensional multifractal level sets
and variational formulae for the entropy in the finite alphabet setting
can be found in \cite{BarreiraSaussolSchmeling:02,BarreiraSaussolSchmeling:02b,BarreiraSaussol:01}.
Therein further interesting examples are provided.

Compairing the finite with the infinite alphabet case we mainly encounter
the following obstacles. The shift space $E_{A}^{\infty}$ is not
even locally compact and hence also $\mathcal{M}\left(E_{A}^{\infty},\sigma\right)$
is not compact with respect to the weak{*}-topology. The function
$Q$ is in general not continuous but only upper semi-continuous with
respect to the weak{*}-topology on $\mathcal{M}\left(E_{A}^{\infty},\sigma\right)$.
This follows from the fact that \begin{equation}
\mu\mapsto\int-I\, d\mu\label{eq:UpperSemi}\end{equation}
 is upper semi-continuous and that $I$ is bounded away from zero.
In Remark \ref{rem:ExampleLSC} we construct as an example measures
$\mu,\mu_{1},\mu_{2},\ldots\in\mathcal{M}\left(E_{A}^{\infty},\sigma\right)$
such that $\mu_{n}\stackrel{*}{\to}\mu$, $\int I\, d\mu<\infty,\int I\, d\mu_{n}<\infty$,
$n\in\N$, but nevertheless $\liminf\int I\, d\mu_{n}>\int I\, d\mu$.
Finally, the entropy map $\mu\mapsto h_{\mu}\left(\sigma\right)$
is not even upper semi-continuous (cf. \cite{JenkinsonMauldinUrbanski:05})
and the pressure function under consideration (cf. Proposition \ref{pro:PressureBetaGibbs})
is only defined in some open region.

We would finally like to remark that the study of higher dimensional
multifractal value sets for infinite GDMS naturally arose in \cite{KesseboehmerStratmann:06}
where a multifractal formalism has been developed in order to study
the fractal geometry of limiting symbols for modular subgroups, which
were introduced by Manin and Marcolli in \cite{ManinMarcolli:02}.

\section{Conformal Graph Directed Markov Systems\label{sec:Conformal-Graph-Directed}}

In this section we begin our study of graph directed Markov systems.
Let us recall the definition of these systems taken from \cite{MauldinUrbanski:03}.
Graph directed Markov systems are based upon a directed multigraph
and an associated incidence matrix, $(V,E,i,t,A)$. The multigraph
consists of a finite set $V$ of vertices and a countable (either
finite or infinite) set of directed edges $E\subset\N$ and two functions
$i,t:E\rightarrow V$. For each edge $e$, $i(e)$ is the initial
vertex of the edge $e$ and $t(e)$ is the terminal vertex of $e$.
The edge goes from $i(e)$ to $t(e)$. Also, a function $A:E\times E\rightarrow\{0,1\}$
is given, called an (edge) incidence matrix. It determines which edges
may follow a given edge. So, the matrix has the property that if $A_{uv}=1$,
then $t(u)=i(v).$ We will consider finite and infinite walks through
the vertex set consistent with the incidence matrix. Thus, we define
the set of infinite admissible words $E_{A}^{\infty}$ on an alphabet
$A$, \[
E_{A}^{\infty}=\{\omega\in E^{\infty}:A_{\omega_{i}\omega_{i+1}}=1\,\,\text{ for all }i\geq1\},\]
 by $E_{A}^{n}$ we denote the set of all subwords of $E_{A}^{\infty}$
of length $n\geq1$, and by $E_{A}^{*}$ we denote the set of all
finite subwords of $E_{A}^{\infty}$. We will consider the left shift
map $\sigma:E_{A}^{\infty}\rightarrow E_{A}^{\infty}$ defined by
$\sigma\left(\omega_{i}\right):=\left(\omega_{i+1}\right)_{i\geq1}$.
Sometimes we also consider this shift as defined on words of finite
length. Given $\omega\in E^{*}$ by $|\omega|$ we denote the length
of the word $\omega$, i.e. the unique $n$ such that $\omega\in E_{A}^{n}$.
If $\omega\in E_{A}^{\infty}$ and $n\geq1$, then \[
\omega|_{n}=\omega_{1}\ldots\omega_{n}.\]
For $\omega\in E_{A}^{\infty}$, or $\omega\in E_{A}^{*}$ with $\left|\omega\right|\geq n$
we will denote with \[
C_{n}\left(\omega\right):=\left\{ x\in E_{A}^{\infty}:x|_{n}=\omega|_{n}\right\} \]
 the \emph{cylinder set of length} $n$ \emph{containing} $\omega$. 

\begin{defn}
\label{def:finiteIrreducible} $E_{A}^{\infty}$ (or equivalently
the GDMS $\Phi$) is called \emph{finitely irreducible} if there exists
a finite set $W\subset E_{A}^{*}$ such that for each $\omega,\eta\in E$
we find $w\in W$ such that the concatenation $\omega w\eta\in E_{A}^{*}$. 
\end{defn}
We recall from the introduction that a \emph{Graph Directed Markov
System} (GDMS) now consists of a directed multigraph and incidence
matrix together with a family of non-empty compact metric spaces $\left(X_{v}\right)_{v\in V}$,
a number $s\in\left(0,1\right)$, and for every $e\in E$, an injective
contraction $\phi_{e}:X_{t(e)}\rightarrow X_{i(e)}$ with a Lipschitz
constant not exceeding $s$. We now describe its limit set. For each
$\omega\in E_{A}^{*}$, say $\omega\in E_{A}^{n}$, we consider the
map coded by $\omega$, \[
\phi_{\omega}:=\phi_{\omega_{1}}\circ\cdots\circ\phi_{\omega_{n}}:X_{t(\omega_{n})}\rightarrow X_{i(\omega_{1})}.\]
 For $\omega\in E_{A}^{\infty}$, the sets $\left\{ \phi_{\omega|_{n}}\left(X_{t(\omega_{n})}\right)\right\} _{n\geq1}$
form a descending sequence of non-empty compact sets and therefore
$\bigcap_{n\geq1}\phi_{\omega|_{n}}\left(X_{t(\omega_{n})}\right)\neq\emptyset$.
Since for every $n\in\N$, $\diam\left(\phi_{\omega|_{n}}\left(X_{t(\omega_{n})}\right)\right)\leq s^{n}\diam\left(X_{t(\omega_{n})}\right)\leq s^{n}\max\{\diam(X_{v}):v\in V\}$,
we conclude that the intersection \[
\bigcap\phi_{\omega|_{n}}\left(X_{t(\omega_{n})}\right)\in X_{i(\omega_{1})}\]
 is a singleton and we denote its only element by $\pi(\omega)$.
In this way we have defined the coding map \begin{equation}
\pi=\pi_{\Phi}:E_{A}^{\infty}\rightarrow X:=\bigoplus_{v\in V}X_{v}\label{eq:pi}\end{equation}
 from $E^{\infty}$ to $\bigoplus_{v\in V}X_{v}$, the disjoint union
of the compact sets $X_{v}$. The set \[
\Lambda=\Lambda_{\Phi}=\pi\left(E_{A}^{\infty}\right)\]
 will be called the \emph{limit set} of the GDMS $\Phi$. 

\begin{defn}
\label{def:Conformal} We call a GDMS \emph{conformal} (CGDMS) if
the following conditions are satisfied.
\begin{itemize}
\item [(a)] For every vertex $v\in V$, $X_{v}$ is a compact connected
subset of a Euclidean space $\R^{D}$ (the dimension $D$ common for
all $v\in V$) and $X_{v}=\overline{\Int(X_{v})}$. 
\item [(b)] \index{Open set condition}\textit{(Open set condition (OSC))}
For all $a,b\in E$, $a\ne b$, \[
\phi_{a}\left(\Int(X_{t(a)}\right)\cap\phi_{b}\left(\Int(X_{t(b)}\right)=\emptyset.\]

\item [(c)] For every vertex $v\in V$ there exists an open connected set
$W_{v}\supset X_{v}$ such that for every $e\in I$ with $t(e)=v$,
the map $\phi_{e}$ extends to a $C^{1}$ conformal diffeomorphism
of $W_{v}$ into $W_{i(e)}$. 
\item [(d)] \index{Cone property}\textit{(Cone property)} There exist
$\gamma,l>0$, $\gamma<\pi/2$, such that for every $x\in X\subset\R^{D}$
there exists an open cone $\Con(x,\gamma,l)\subset\Int(X)$ with vertex
$x$, central angle of measure $\gamma$, and altitude $l$. 
\item [(e)] There are two constants $L\geq1$ and $\alpha>0$ such that
\[
\left||\phi_{e}'(y)|-|\phi_{e}'(x)|\right|\leq L\Vert(\phi_{e}')^{-1}\Vert^{-1}\Vert y-x\Vert^{\alpha}\]
 for every $e\in I$ and every pair of points $x,y\in X_{t(e)}$,
where $|\phi_{\omega}'(x)|$ means the norm of the derivative. 
\end{itemize}
\end{defn}
The following remarkable fact was proved in \cite{MauldinUrbanski:03}.

\begin{prop}
\label{p1.033101} If $D\geq2$ and a family $\Phi=\left(\phi_{e}\right)_{e\in I}$
satisfies conditions (a) and (c), then it also satisfies condition
(e) with $\alpha=1$. 
\end{prop}
The following rather straightforward consequence of (e) was proved
in \cite{MauldinUrbanski:03}.

\begin{lem}
\label{l2.033101} If $\Phi=\left(\phi_{e}\right)_{e\in I}$ is a
CGDMS, then for all $\omega\in E^{*}$ and all $x,y\in W_{t(\omega)}$,
we have \[
\left|\log|\phi_{\omega}'(y)|-\log|\phi_{\omega}'(x)|\right|\leq\frac{L}{1-s}\Vert y-x\Vert^{\alpha}.\]

\end{lem}
As a straightforward consequence of (e) we get the following.

\begin{itemize}
\item [(f)] (Bounded distortion property). There exists $K\geq1$ such
that for all $\omega\in E^{*}$ and all $x,y\in X_{t(\omega)}$ \[
|\phi_{\omega}'(y)|\leq K|\phi_{\omega}'(x)|.\]

\end{itemize}
Next we define the geometrical potential function associated with
$\Phi$ by \[
I=I_{\Phi}:E_{A}^{\infty}\rightarrow\R^{+},\: I_{\Phi}\left(\omega\right):=-\log\left|\phi'_{\omega_{1}}\left(\pi\left(\sigma\left(\omega\right)\right)\right)\right|.\]

It was proved in \cite{MauldinUrbanski:03} that for each $t\geq0$
the following limit exists (possibly be equal to $+\infty$). \[
\mathfrak{p}\left(t\right):=\P(-tI):=\lim_{n\rightarrow\infty}\frac{1}{n}\log\sum_{\omega\in E_{A}^{n}}\exp\left(-tS_{n}I\left(\omega\right)\right)=\lim_{n\rightarrow\infty}\frac{1}{n}\log\sum_{\omega\in E_{A}^{n}}||\phi_{\omega}'||^{t},\]

where $S_{n}I\left(\omega\right):=\sup_{\left\{ x:x|_{n}=\omega\right\} }\sum_{i=0}^{n-1}I\left(\sigma^{i}x\right)$.
This number is called the \emph{topological pressure} of the parameter
$t$. The function $\mathfrak{p}$ is always non-increasing and convex.
In \cite{MauldinUrbanski:03} a useful parameter associate with a
CGDMS has been introduced. Namely, \[
\theta(\Phi):=\inf\{ t:\mathfrak{p}(t)<\infty\}=\sup\{ t:\mathfrak{p}(t)=\infty\}.\]
 Let $\Fin(E)$ denote the set of all finite subsets of $E$ and for
$F\in\Fin\left(E\right)$ we define the subsystem $\Phi_{F}:=\left(\phi_{i}\right)_{i\in F}$
of $\Phi$. The following characterization of $h_{\Phi}=\HD\left(\Lambda_{\Phi}\right)$,
the Hausdorff dimension of the limit set $\Lambda_{\Phi}$, being
a variant of Bowen's formula, was proved in \cite{MauldinUrbanski:03}
as Theorem~4.2.13.

\begin{thm}
\label{gdms4213} If the CGDMS $\Phi$ is finitely irreducible, then
\[
\HD\left(\Lambda_{\Phi}\right)=\inf\left\{ t\geq0:\mathfrak{p}(t)<0\right\} =\sup\left\{ h_{\Phi_{F}}:F\in\Fin(I)\right\} \geq\theta(\Phi).\]
 If $\mathfrak{p}(t)=0$, then $t$ is the only zero of the function
$\mathfrak{p}(t)$, $t=\HD(\Lambda_{\Phi})$ and the system $\Phi$
is called \emph{regular}. 
\end{thm}
In fact it was assumed in \cite{MauldinUrbanski:03} that the system
$\Phi$ is finitely primitive but the proof can be easily improved
to this slightly more general setting. It will be convenient for us
to recall and make use of the following definitions.

\begin{defn}
\label{d4.2.1.} A CGDMS is said to be \emph{strongly regular} if
there exists $t\geq0$ such that $0<\mathfrak{p}\left(t\right)<\infty$.
A family $\left(\phi_{i}\right)_{i\in F}$ is said to be a \emph{co-finite
subsystem} of a system of $\Phi=\left(\phi_{i}\right)_{i\in E}$ if
$F\subset E$ and the difference $E\setminus F$ is finite. 
\begin{defn}
\label{d4.2.2} A CGDMS is said to be \emph{co-finitely regular} if
each of its co-finite subsystem is regular. 
\end{defn}
\end{defn}
The following fact relating all these three notions can be found in
\cite{MauldinUrbanski:03}.

\begin{prop}
\label{p1121005} Each co-finitely regular system is strongly regular
and each strongly regular system is regular. 
\end{prop}
Note that the system $\Phi$ is strongly regular if and only if $\HD(\Lambda_{\Phi})>\theta(\Phi)$.

\section{Higher-dimensional thermodynamic formalism}

From now on let us assume that the CGDMS is infinite and co-finitely
regular, i.e. $\card\left(E\right)=\infty$ and hence $\theta\left(\Phi\right)\geq0$.
Recall that for $\omega,\tau\in E_{A}^{\infty}$, we define $\omega\wedge\tau\in E_{A}^{\infty}\cup E_{A}^{*}$
to be the longest initial block common to both $\omega$ and $\tau$.
We say that a function $f:E_{A}^{\infty}\rightarrow\R$ is \emph{H\"{o}lder
continuous} with an exponent $\alpha>0$ if \[
v_{\alpha}(f):=\sup\left\{ V_{\alpha,n}(f):n\geq1\right\} <\infty,\]
where\[
V_{\alpha,n}(f)=\sup\left\{ |f(\omega)-f(\tau)|\e^{\alpha(n-1)}:\omega,\tau\in E_{A}^{\infty}\text{ and }|\omega\wedge\tau|\geq n\right\} .\]
 For every $\alpha>0$ let $\K_{\alpha}$ be the set of all real-valued
H\"{o}lder continuous (not necessarily bounded) functions on $E_{A}^{\infty}$.
Set \[
\K_{\alpha}^{s}:=\left\{ f\in\K_{\alpha}:\sum_{e\in E}\exp\left(\sup\left(f|_{C_{1}\left(e\right)}\right)\right)<+\infty\right\} .\]
 Each member of $\K_{\alpha}^{s}$ is called an $\alpha$\emph{-H\"{o}lder
summable potential}.

For fixed $d\in\N$ let $J:E_{A}^{\infty}\rightarrow\R^{d}$ such
that $J_{i}\in\K_{\alpha}$ is a \emph{bounded} Hölder continuous
function for $i=1,\ldots,d$. The following proposition will be of
central importance throughout this paper.

\begin{prop}
\label{pro:PressureBetaGibbs}Each member of the family $\left(\left\langle t,J\right\rangle -\beta I:t\in\R^{d},\beta>\theta\right)$
is an element of $\K_{\alpha}^{s}$. The pressure functional \[
p:\R^{d}\times\left(\theta,\infty\right)\rightarrow\R,\; p\left(t,\beta\right):=\P\left(\left\langle t,J\right\rangle -\beta I\right)\]
 is a well-defined, real-analytic, convex function. For each $t\in\R^{d}$
there exists a unique number $\beta\left(t\right)$ such that $p\left(t,\beta\left(t\right)\right)=0$.
Also $t\mapsto\beta\left(t\right)$ defines a real-valued, real-analytic
convex function on $\R^{d}$. Its gradient is given by \begin{equation}
\nabla\beta\left(t\right)=\frac{1}{\int I\, d\mu_{t}}\int J\, d\mu_{t},\label{eq:gradient}\end{equation}
 where $\mu_{t}=\mu_{t,\beta\left(t\right)}$ denotes the unique invariant
Gibbs measure for the potential $\left\langle t,J\right\rangle -\beta\left(t\right)I$,
i.e. there exists $C>0$ such that for all $\omega\in E_{A}^{\infty}$
we have\begin{equation}
C^{-1}\leq\frac{\mu_{t}C_{n}\left(\omega\right)}{\exp S_{n}\left(\left\langle t,J\right\rangle -\beta\left(t\right)I\right)\left(\omega\right)}\leq C.\label{eq:Gibbs}\end{equation}
 
\end{prop}
\begin{proof}
The properties of the family $\left(\left\langle t,J\right\rangle -\beta I:t\in\R^{d},\beta>\theta\right)$
follows immediately from the boundedness of $J$ and the Hölder continuity
of $I$. From \cite{MauldinUrbanski:03} we then know that $p$ is
a well-defined and real-analytic function. Since the system is infinite
and co-finitely regular we have $\lim_{\beta\searrow\theta}p\left(t,\beta\right)=\infty$.
Furthermore, for every $t\in\R^{d}$, we have \[
\partial_{\beta}p\left(t,\beta\right)=-\int I\, d\mu_{t,\beta}\leq\log s<0.\]
 Hence, for every $t\in\R^{d}$, $\beta\mapsto p\left(t,\beta\right)$
is a strictly decreasing function and $\lim_{\beta\rightarrow+\infty}p\left(t,\beta\right)=-\infty$.
From this we conclude that for each $t\in\R^{d}$ there exists a unique
number $\beta\left(t\right)>\theta$ such that $p\left(t,\beta\left(t\right)\right)=0$.
By the implicit function theorem $\beta:\R^{d}\rightarrow\R$ is real-analytic
and convex. The formula for the gradient of $\beta$ follows again
from the implicit function theorem. 
\end{proof}
\begin{lem}
\label{lem:tight} Any set of measures $M\subset\mathcal{M}\left(E_{A}^{\infty},\sigma\right)$
such that $\sup_{\mu\in M}\int I\, d\mu<\infty$ is tight. 
\end{lem}
\begin{proof}
For every $i,\ell\in\N$ put \[
E_{i,\ell}:=\left\{ \omega\in E_{A}^{\infty}:\omega_{i}\geq\ell\right\} .\]
Then we have for all $\mu\in M$\[
\textrm{const.}\geq\int I\, d\mu\geq\int_{E_{1,\ell}}I\, d\mu\geq\mu\left(E_{1,\ell}\right)\inf_{E_{1,\ell}}I.\]
Combining this with the fact that $E_{i,\ell}\subset\sigma^{-i+1}\left(E_{1,\ell}\right)$
and that $\mu$ is $\sigma$--invariant we get \[
\mu\left(E_{i,\ell}\right)\leq\mu\left(E_{1,\ell}\right)\leq\frac{\textrm{const.}}{\inf_{E_{1,\ell}}I}.\]
Now, fix $\epsilon>0$. It follows from the above estimate that for
every $i\geq1$ there exists $\ell_{i}\geq1$ such that $\mu\left(E_{i,\ell_{i}}\right)<2^{-i}\epsilon$.
Then $\bigcup_{i=1}^{\infty}E_{i,\ell_{i}}$ is the complement of
the compact set $\left\{ \omega\in E_{A}^{\infty}:\forall i\in\N:\:\omega_{i}<\ell_{i}\right\} $
and \[
\mu\left(\bigcup_{i=1}^{\infty}E_{i,\ell_{i}}\right)\leq\sum_{i=1}^{\infty}\mu\left(E_{i,\ell_{i}}\right)\leq\sum_{i=1}^{\infty}2^{-i}\epsilon=\epsilon\quad\textrm{for all}\quad\mu\in M.\]
From this the tightness of $M$ follows.
\end{proof}

\section{Proofs}

Clearly, since $I>-\log s$ and $J$ is bounded, the sets $K,L,M\subset\R^{d}$
are all bounded.

Let us first show that \begin{equation}
0\in K\cap L.\label{eq:0inKandL}\end{equation}
 To see this we construct a Bernoulli measure $\mu_{p}$ (which is
invariant and ergodic) with probability vector $p:=\left(p_{i}\right)$
chosen in such a way that $\sum p_{i}\inf I|_{C_{1}\left(i\right)}=\infty$,
which is always possible. Then $\left|Q\left(\mu_{p}\right)\right|=\left|\frac{\int J\, d\mu_{p}}{\int I\, d\mu_{p}}\right|\leq$`$\frac{\textrm{const.}}{\infty}$'$=0$
and for $\mu_{p}$--almost all points $\omega\in E_{A}^{\infty}$
we have by the ergodic theorem that $\lim_{n\rightarrow\infty}\frac{S_{n}J\left(\omega\right)}{S_{n}I\left(\omega\right)}=Q\left(\mu_{p}\right)=0$.

\subsection{Proof of Theorem \ref{thm:Main1} \textmd{}}

~

The inclusion $M\subset L$ follows immediately from the definitions
and (\ref{eq:gradient}).

\subsubsection{Proof of `$\Int L\subset M$'.}

~

Here we follow some ideas from \cite{BarreiraSaussolSchmeling:02}.
Let $\alpha\in\Int L$. Then there exists $r>0$ such that $B_{r}\left(\alpha\right)\subset L$.
Let \[
D_{\alpha}:=\sup_{\mu\in Q^{-1}\left(\alpha\right)}\left\{ \frac{h\left(\mu\right)}{\int I\, d\mu}\right\} ,\,\alpha\ne0.\]
By the variational principle (cf. \cite[Theorem 2.1.7]{MauldinUrbanski:03})
it follows that \[
0=\mathcal{P}\left(\HD\left(\Lambda_{\Phi}\right)I\right)\geq h\left(\mu\right)-\HD\left(\Lambda_{\Phi}\right)\int I\, d\mu\]
 and hence $D_{\alpha}$ is always dominated by the Hausdorff dimension
$\HD\left(\Lambda_{\Phi}\right)$. Let us now consider the family
of potentials \[
\left(\left\langle q,J\right\rangle -\left(\left\langle q,\alpha\right\rangle +D_{\alpha}\right)I:q\in G_{\alpha}\right),\]
 where $G_{\alpha}:=\left\{ x\in\R^{d}:\left\langle x,\alpha\right\rangle +D_{\alpha}>\theta\right\} $.
Firstly, we show for all $q\in G_{\alpha}$ \[
p_{\alpha}\left(q\right):=P\left(\left\langle q,J\right\rangle -\left(\left\langle q,\alpha\right\rangle +D_{\alpha}\right)I\right)\geq0.\]
 Indeed, by the variational principle we have\begin{eqnarray*}
p_{\alpha}\left(q\right) & \geq & \sup_{\mu\in Q^{-1}\left(\alpha\right)}\left\{ h_{\mu}+\left\langle q,\int J\, d\mu\right\rangle -\left(\left\langle q,\alpha\right\rangle +D_{\alpha}\right)\int I\, d\mu\right\} \\
 & \geq & \sup_{\mu\in Q^{-1}\left(\alpha\right)}\left\{ \int I\, d\mu\left(\frac{h_{\mu}}{\int I\, d\mu}-D_{\alpha}\right)\right\} \\
 & \gg & \sup_{\mu\in Q^{-1}\left(\alpha\right)}\left\{ \frac{h_{\mu}}{\int I\, d\mu}-D_{\alpha}\right\} =0.\end{eqnarray*}
 Here we used the fact that $\sup_{\mu\in Q^{-1}\left(\alpha\right)}\int I\, d\mu$
is bounded above by some constant for $\alpha\ne0$.

Next we show that the infimum $\inf_{q\in G_{\alpha}}p_{\alpha}\left(q\right)$
is attained at some point $q\in G_{\alpha}$. This follows from the
fact that $p_{\alpha}\left(q_{n}\right)$ diverges to infinity whenever
either $\left\langle q_{n},\alpha\right\rangle +D_{\alpha}\rightarrow\theta$,
$n\rightarrow\infty$, which is clear, or $\left|q_{n}\right|\rightarrow\infty$
which can be seen as follows. For $q=\left(q_{1},\ldots,q_{d}\right)$
let $\beta_{i}:=\alpha_{i}+r/2\sgn q_{i}$, $i=1,\ldots,d$ and $\mu\in Q^{-1}\left(\beta\right)$.
Then we have \begin{eqnarray}
p_{\alpha}\left(q\right) & \geq & h_{\mu}+\left\langle q,\int\left(J-\alpha I\right)\, d\mu\right\rangle -D_{\alpha}\int I\, d\mu\nonumber \\
 & = & h_{\mu}+\left\langle q,\left(\beta-\alpha\right)\right\rangle \int I\, d\mu-D_{\alpha}\int I\, d\mu\nonumber \\
 & = & h_{\mu}+\left(\left\langle q,\left(\beta-\alpha\right)\right\rangle -D_{\alpha}\right)\int I\, d\mu\nonumber \\
 & \geq & \left(\frac{r}{2}\sum_{i=1}^{d}\left|q_{i}\right|-\HD\left(\Lambda_{\Phi}\right)\right)\left(-\log s\right).\label{eq:lowerbound}\end{eqnarray}
 Now the right hand side diverges to infinity for $\left|q\right|\rightarrow\infty$
showing that the infimum must be attained in some uniformly bounded
region, say in $q_{\alpha}\in G_{\alpha}$. Since $p_{\alpha}$ is
real-analytic on $G_{\alpha}$ we have\begin{equation}
0=\nabla p_{\alpha}\left(q_{\alpha}\right)=\int J-\alpha I\, d\mu_{q_{\alpha}}\implies\alpha=\frac{\int J\, d\mu_{q_{\alpha}}}{\int I\, d\mu_{q_{\alpha}}},\label{eq:derivertive}\end{equation}
where $\mu_{q_{\alpha}}$ is the Gibbs measure for the potential $\left(\left\langle q_{\alpha},J\right\rangle -\left(\left\langle q_{\alpha},\alpha\right\rangle +D_{\alpha}\right)I\right)$.
Hence, by the variational principle, we have\[
p\left(q_{\alpha}\right)=h_{\mu_{q_{\alpha}}}-D_{\alpha}\int I\, d\mu_{q_{\alpha}}\geq0\implies D_{\alpha}\leq\frac{h_{\mu_{q_{\alpha}}}}{\int I\, d\mu_{q_{\alpha}}}.\]
 By definition of $D_{\alpha}$ we actually have equality and hence
$p\left(q_{\alpha}\right)=0$. This, (\ref{eq:derivertive}) and Proposition
\ref{pro:PressureBetaGibbs} show that $\alpha=\nabla\beta\left(q_{\alpha}\right)$
for $\alpha\ne0$.

Now we consider the case $\alpha=0$. If $0\in\partial L$ then nothing
has to be shown. If $0\in\Int\left(L\right)$ then also $B_{r}\left(0\right)\subset\Int\left(L\right)$
for some small $r>0$ and by the above we have $B_{r}\left(0\right)\setminus\left\{ 0\right\} \subset M$.
The convexity of $\beta$ then implies that $\beta$ has a minimum
in $\mathbb{R}^{d}$ and since $\beta$ is real-analytic this minimum
is unique, say in $t_{0}$, with $\nabla\beta\left(t_{0}\right)=0$.
This shows that $0\in M$.

\subsubsection{Proof of `$K\subset\overline{L}$'.}

~

Let $\alpha\in K$. In view of (\ref{eq:0inKandL}) we may assume
without loss of generality that $\alpha$ is not equal to zero. Then
there exist $\omega\in E_{A}^{\infty}$ such that $\lim_{n\rightarrow\infty}\frac{S_{n}J\left(\omega\right)}{S_{n}I\left(\omega\right)}=\alpha$.
Now, for every $n\in\N$ there exists by the finite irreducibility
condition a word $w_{n}\in W$ such that the periodic element $x_{n}:=\left(\omega|_{n}w_{n}\right)^{\infty}$
belongs to $E_{A}^{\infty}$. The point $x_{n}$ gives then rise to
the invariant probability measure \[
\mu_{n}:=\frac{1}{k_{n}}\sum_{j=0}^{k_{n}-1}\delta_{\sigma^{j}x_{n}}\in\mathcal{M}\left(E_{A}^{\infty},\sigma\right),\]
 where $k_{n}:=n+\left|w_{n}\right|.$ By the Hölder continuity (bounded
distortion) and the finiteness of $W$ we estimate\begin{eqnarray}
\frac{\int J\, d\mu_{n}}{\int I\, d\mu_{n}} & = & \frac{S_{k_{n}}J\left(x_{n}\right)}{S_{k_{n}}I\left(x_{n}\right)}=\frac{S_{n}J\left(x_{n}\right)+\textrm{O}\left(1\right)}{S_{n}I\left(x_{n}\right)+\textrm{O}\left(1\right)}=\frac{S_{n}J\left(\omega\right)+\textrm{O}\left(1\right)}{S_{n}I\left(\omega\right)+\textrm{O}\left(1\right)},\label{eq:ApproxPeriod}\end{eqnarray}
where $\textrm{O}$ denotes the corresponding Landau symbol. Since
$S_{n}I$ growth at least like $-n\log s$ the above quotient converges
to $\alpha$ as $n\to\infty$. This shows $\alpha\in\overline{L}$.

\subsubsection{Proof of compactness of $K$}

~

Since we know that $K$ is bounded, we are left to show that $K$
is closed. By (\ref{eq:0inKandL}) we have $0\in K$. Hence, without
loss of generality we may consider a sequence $\left(\alpha_{k}\right)\in K^{\N}$
converging to $\alpha\in\R^{d}\setminus\left\{ 0\right\} $. We are
going to construct inductively an element $\omega\in E_{A}^{\infty}$
such that $\lim\frac{S_{n}J\left(\omega\right)}{S_{n}I\left(\omega\right)}=\alpha$.
Fix a sequence $\epsilon_{k}\searrow0$ such that $\left|\alpha_{k}-\alpha\right|<\epsilon_{k}/2$.
Using the observation in (\ref{eq:ApproxPeriod}) we find for each
$k\in\N$ a periodic element $x_{k}=p_{k}^{\infty}\in E_{A}^{\infty}$,
$m_{k}:=\left|p_{k}\right|$, such that $\left|\frac{S_{m_{k}}J\left(x_{k}\right)}{S_{m_{k}}I\left(x_{k}\right)}-\alpha_{k}\right|<\frac{\epsilon_{k}}{2}$
which gives\[
\left|\frac{S_{m_{k}}J\left(x_{k}\right)}{S_{m_{k}}I\left(x_{k}\right)}-\alpha\right|<\epsilon_{k}\]
 We begin the induction by defining $\omega_{1}:=p_{1}^{l_{1}}w_{1}$
with $l_{1}=1$.

Suppose we have already defined $\omega_{k}:=p_{1}^{l_{1}}w_{1}\cdots p_{k}^{l_{k}}w_{k}$
for some $k\in\N$ and let \[
N_{k}:=\sum_{i=1}^{k}l_{i}m_{i}+\left|w_{i}\right|.\]
 Then choose $w_{k+1}\in W$ such that $\omega_{k}w_{k+1}p_{k+1}\in E_{A}^{*}$
and $l_{k+1}\in\N$ large enough such that \[
\frac{1}{-\log s}\cdot\frac{1}{l_{k+1}k_{k+1}}\max\left\{ S_{m_{k+2}}I\left(x_{k+2}\right),S_{m_{k+2}}\left|J\left(x_{k+2}\right)\right|\right\} \leq\epsilon_{k+1}.\]
 In this way we define inductively the infinite word $\omega:=\left(p_{i}^{l_{i}}w_{i}\right)_{i=1}^{\infty}\in E_{A}^{\infty}.$ 

We will need the following observation. Suppose we have two sequences
$\left(a_{n}\right)\in\left(\R^{d}\right)^{\N}$ and $\left(b_{n}\right)\in\left(\R^{+}\right)^{\N}$
such that $b_{n}^{-1}a_{n}\to\alpha$ and $\liminf_{n}b_{n}>0$. Define
$A_{N}:=\sum_{k=1}^{N}a_{k}$ and $B_{N}:=\sum_{k=1}^{N}b_{k}$. Let
$\textrm{o}$ denote the corresponding Landau symbol. Then for any
two sequences $\left(c_{n}\right)$ and $\left(d_{n}\right)$ given
by $c_{n}:=A_{k_{n}}+\textrm{o}\left(B_{k_{n}}\right)$ and $d_{n}:=B_{k_{n}}+\textrm{o}\left(B_{k_{n}}\right)$
for some sequence $\left(k_{n}\right)\in\N^{\N}$ tending to infinity,
we have $d_{n}^{-1}c_{n}\to\alpha$. 

For $n\in\N$ we define a sequence $\left(k_{n}\right)\in\N^{\N}$
such that $N_{k_{n}}\leq n<N_{k_{n}+1}$, and $r_{n}$, $\ell_{n}$
such that $n=N_{k_{n}}+q_{n}\cdot m_{k_{n}+1}+r_{n}$ with $0\leq r_{n}\leq m_{k_{n}+1}$
and $0\leq q_{n}\leq\ell_{k_{n}+1}$. Then applying the above observation
to \[
S_{n}J\left(\omega\right)=\sum_{i=1}^{k_{n}}l_{i}S_{m_{i}}J\left(x_{i}\right)+q_{n}S_{m_{k_{n}+1}}J\left(x_{k_{n}+1}\right)+\textrm{O}\left(S_{r_{n}}J\left(x_{k_{n}+1}\right)\right)+\textrm{O}\left(k_{n}\right)\]
 and \[
S_{n}I\left(\omega\right)=\sum_{i=1}^{k_{n}}l_{i}S_{m_{i}}I\left(x_{i}\right)+q_{n}S_{m_{k_{n}+1}}I_{n}\left(x_{k_{n}+1}\right)+\textrm{O}\left(S_{r_{n}}I\left(x_{k_{n}+1}\right)\right)+\textrm{O}\left(k_{n}\right)\]
 and by observing the definition of $\left(\ell_{k}\right)$ the claim
follows.

\subsubsection{Proof of `$\overline{M}\subset K$'}

~

Since we have already seen that $K$ is closed it suffices to prove
$M\subset K$. Let $\alpha=\nabla\beta\left(t\right)\in M$. Then
for the Gibbs measure $\mu_{t}$ we have $\alpha=\left(\int I\, d\mu_{t}\right)^{-1}\int J\, d\mu_{t}$.
By the ergodicity of $\mu_{t}$ we have for $\mu_{t}$-a.e. $x\in E_{A}^{\infty}$
that $\lim_{n\to\infty}\frac{S_{n}J\left(x\right)}{S_{n}I\left(x\right)}=\alpha$
and hence $\alpha\in K$.

\subsection{Proof of Theorem \ref{thm:Main2} \textmd{}}

~

For $J_{i}$ linearly independent as cohomology classes it is well
known that $\beta$ is strictly convex, or equivalently the Hessian
$\textrm{Hess}\left(\beta\right)$ is strictly positive definite.
From this it follows that $\nabla\beta:\R^{d}\rightarrow\nabla\beta\left(\R^{d}\right)$
is a diffeomorphism and hence $M:=\nabla\beta\left(\R^{d}\right)$
is an open and connected subset of $\R^{d}$.

\subsubsection{Proof of `$L\subset\overline{\Int L}$'}

~

Since $\emptyset\ne M\subset\Int L$ we may use similar arguments
as in \cite{BarreiraSaussolSchmeling:02} to prove that in this situation
we have $L\subset\overline{\Int L}$. Let $\alpha=Q\left(m_{0}\right)\in L$.
Since $\Int L$ is not empty we find measures $m_{1},\ldots,m_{d}$
such that $\left(\frac{\int J\, dm_{1}}{\int I\, dm_{1}}-\alpha,\ldots,\frac{\int J\, dm_{d}}{\int I\, dm_{d}}-\alpha\right)$
form a basis of $\R^{d}$. For $p\in\Delta:=\left\{ u\in\R^{d}:u_{i}\geq0\,\wedge\,\sum_{l}u_{l}\leq1\right\} $
we define $\mu_{p}:=\sum_{l=1}^{d}p_{l}m_{l}+\left(1-\sum p_{l}\right)m_{0}$.
Then the derivative of $b:\Delta\to\R^{d}$, $b\left(p\right)=Q\left(\mu_{p}\right)$
is given by \begin{eqnarray*}
\frac{\partial b_{i}}{\partial p_{j}}\left(0\right) & = & \frac{\left(\int J_{i}\, dm_{j}-\int J_{i}\, dm_{0}\right)}{\int I\, dm_{0}}-\frac{\int J_{i}\, dm_{0}\left(\int I\, dm_{j}-\int I\, dm_{0}\right)}{\left(\int I\, dm_{0}\right)^{2}}\\
 & = & \frac{\int J_{i}\, dm_{j}-\alpha_{i}\int I\, dm_{j}}{\int I\, dm_{0}}.\end{eqnarray*}
By our assumption $\left(\frac{\int I\, dm_{1}}{\int I\, dm_{0}}\left(\frac{\int J\, dm_{1}}{\int I\, dm_{1}}-\alpha\right),\ldots,\frac{\int I\, dm_{d}}{\int I\, dm_{0}}\left(\frac{\int J\, dm_{d}}{\int I\, dm_{d}}-\alpha\right)\right)$
are linearly independent, and hence $\frac{db}{dp}$ is invertible.
This shows that there exists an open set $U\subset\Delta$ such that
$0\in\overline{U}$ and $b:U\to b\left(U\right)$ is a diffeomorphism.
We finish the argument by observing that $\alpha\in\overline{b\left(U\right)}\subset\overline{\Int L}$.

\subsubsection{Proof of `$L=\overline{M}$ for $0\in M$' \label{sub:Proof-of-L=3DoverbarM}}

~

Since $\Int L\subset M\subset L\subset\overline{\Int L}$ we have
$\overline{M}=\overline{L}$. It now suffices to show that $\overline{M}\subset L$
since the inclusion $\overline{L}=\overline{M}\subset L\subset\overline{L}$
would then imply that $L=\overline{M}$. To see that indeed $\overline{M}\subset L$
we proceed as follows. Recall that by our assumptions $0\in M$. Let
$\alpha\in\partial M$, which is then necessary different from $0$.
Let $\left(\alpha_{n}\right)\in M^{\N}$ be a sequence converging
to $\alpha$. For this sequence we find a sequence of Gibbs measures
$\left(\mu_{s_{k}}\right)$ (for the potential $\left\langle s_{k},J\right\rangle -\beta\left(s_{k}\right)I$)
such that $\left(\int I\, d\mu_{s_{k}}\right)^{-1}\int J\, d\mu_{s_{k}}$
converges to $\alpha$ and $\left|s_{k}\right|\rightarrow\infty$.
In particular, we have that $\left(\int I\, d\mu_{s_{k}}\right)$
is bounded . By Lemma \ref{lem:tight} this sequence of measures is
tight and hence there is a weak convergent subsequence $\mu_{t_{k}}\rightarrow\mu\in\mathcal{M}\left(E_{A}^{\infty},\sigma\right)$.
Now we have to show that $Q\left(\mu\right)=\alpha$. We clearly have
$\int J\, d\mu_{t_{k}}\rightarrow\int J\, d\mu=v$ since $J$ is bounded
and by (\ref{eq:UpperSemi}) we have\begin{equation}
\liminf_{k\rightarrow\infty}\int I\, d\mu_{t_{k}}\geq\int I\, d\mu.\label{eq:liminfGEQ}\end{equation}
 To show furthermore that $\limsup_{k\rightarrow\infty}\int I\, d\mu_{t_{k}}\leq\int I\, d\mu$
we make use of the variational principle. Indeed we have\[
0=h\left(\mu_{t_{k}}\right)+\left\langle t_{k},\int J\, d\mu_{t_{k}}\right\rangle -\beta\left(t_{k}\right)\int I\, d\mu_{t_{k}}\geq h_{\mu}+\left\langle t_{k},\int J\, d\mu\right\rangle -\beta\left(t_{k}\right)\int I\, d\mu\]
 This gives \begin{equation}
\int I\, d\mu\geq\frac{h_{\mu}-h_{\mu_{t_{k}}}}{\beta\left(t_{k}\right)}+\left\langle \beta\left(t_{k}\right)^{-1}\cdot t_{k},\int J\, d\mu-\int J\, d\mu_{t_{k}}\right\rangle +\int I\, d\mu_{t_{k}}.\label{eq:lowerBoundIdmu}\end{equation}
Since\[
0\leq h_{\mu_{t}}\leq\HD\left(\Lambda_{\Phi}\right)\int I\, d\mu_{t}\]
 and since $\int I\, d\mu_{t_{k}}$ is bounded we conclude that $h_{\mu_{t}}$
is bounded. The assumption $0\in M$ implies that $\left|t_{k}\right|/\beta\left(t_{k}\right)$
is bounded (cf. \cite[Theorem 3.26]{RochafellarWets:98}). As $J$
is a bounded function we have $\int J\, d\mu-\int J\, d\mu_{t_{k}}\rightarrow0$
and hence taking limits in (\ref{eq:lowerBoundIdmu}) gives $\limsup_{k\rightarrow\infty}\int I\, d\mu_{t_{k}}\leq\int I\, d\mu$.
Combining this with (\ref{eq:liminfGEQ}) finally proves $Q\left(\mu\right)=\alpha$.

The remaining parts of the theorem are an immediate application of
Theorem \ref{thm:Main1}. 

\begin{rem}
\label{rem:ExampleLSC} We would like to emphasis that, unlike above,
for $\mu,\mu_{1},\mu_{2},\ldots\in\mathcal{M}\left(\Sigma_{E}^{\infty},\sigma\right)$
with $\mu_{n}\stackrel{*}{\to}\mu$ and $\int I\, d\mu<+\infty$,
$\int I\, d\mu_{n}<\infty,$ $n\in\N$, the convergence $\int I\, d\mu_{n}\to\int I\, d\mu$
does not hold in general. For a counter example we consider again
the iterated function system $\Phi$ generated by continued fractions
as given in (\ref{eq:CF}). Fix $M>0$ to be large; $c_{n}:=1-M/\log\left(n\right)$,
$n>\exp\left(M\right)$; $S:=\sum_{k=1}^{\infty}k^{-2}$ and $\left(p_{k}^{\left(n\right)}\right)_{k\in\N}$
a probability vector such that\[
p_{k}^{\left(n\right)}:=\left\{ \begin{array}{lll}
0 & \textrm{ if } & k>n,\\
S^{-1}c_{n}k^{-2} & \textrm{ if } & k<n,\\
\left(1-S^{-1}c_{n}\sum_{j=1}^{n-1}j^{-2}\right) & \textrm{ if } & k=n.\end{array}\right.\]
Let $\mu_{n}$ be the Bernoulli measure associated with $\left(p_{k}^{\left(n\right)}\right)_{k\in\N}$.
Then $\left(\mu_{n}\right)_{n>\exp M}$ converges weakly to the Bernoulli
measure associated with the probability vector $\left(S^{-1}k^{-2}\right)_{k\in\N}$.
We then have on the one hand\begin{eqnarray*}
\int I\, d\mu & \leq & 2S^{-1}\sum_{k\in\N}\frac{\log\left(k+1\right)}{k^{2}}<\infty\\
\int I\, d\mu_{n} & \leq & 2\left(1-S^{-1}c_{n}\sum_{j=1}^{n-1}j^{-2}\right)\log\left(n+1\right)+\sum_{k=1}^{n-1}2S^{-1}c_{n}\frac{\log\left(k+1\right)}{k^{2}}\\
 & \leq & 2\left(1-S^{-1}c_{n}\left(S-\frac{c}{n}\right)\right)\log\left(n+1\right)+2S^{-1}\sum_{k=1}^{n-1}c_{n}\frac{\log\left(k+1\right)}{k^{2}}\\
 & \leq & 2\left(\frac{M}{\log n}+\frac{cS^{-1}}{n}\right)\log\left(n+1\right)+2S^{-1}\sum_{k=1}^{n-1}c_{n}\frac{\log\left(k+1\right)}{k^{2}}\\
 & \leq & 4M+2S^{-1}\sum_{k=1}^{n-1}c_{n}\frac{\log\left(k+1\right)}{k^{2}}<\infty,\end{eqnarray*}
for some constant $c>0$ and for $M$ sufficiently large. On the other
hand we have\begin{eqnarray*}
\int I\, d\mu_{n} & \geq & \int_{C_{1}\left(\left[1,\ldots\right]\right)}I\, d\mu\geq2\left(1-S^{-1}c_{n}\sum_{j=1}^{\infty}j^{-2}\right)\log n\\
 & \geq & 2\left(1-c_{n}\right)\log n=2M\geq2\int I\, d\mu,\end{eqnarray*}
for $M$ large enough. Hence in this example, we have \[
\liminf_{k\to\infty}\int I\, d\mu_{n}>\int I\, d\mu.\]

\end{rem}
For the proof of Theorem \ref{thm:MF} we need the following lemma
from convex analysis adapted to our situation. For an extended real-valued
function $f:\R^{d}\to\overline{\R}$ we define the effective domain\[
\dom f:=\left\{ x\in\R^{d}:f\left(x\right)>-\infty\right\} .\]

\begin{lem}
Under the conditions of Theorem \ref{thm:MF} we have $M=\Int\left(\dom\hat{\beta}\right)$
is a non-empty convex set. Furthermore, for each $a\in M$ and $\alpha\in\partial M$,
we have \begin{equation}
\lambda\alpha+\left(1-\lambda\right)a\in M\:\textrm{for all }\lambda\in[0,1),\:\textrm{and }\,\lim_{\lambda\to1}\hat{\beta}\left(\lambda\alpha+\left(1-\lambda\right)a\right)=\hat{\beta}\left(\alpha\right).\label{eq:Continuity}\end{equation}
In particular, we have \begin{equation}
\overline{M}=\dom\hat{\beta}.\label{eq:Mdom}\end{equation}
 
\end{lem}
\begin{proof}
Since $\Int\left(\dom\hat{\beta}\right)\subset M\subset\dom\hat{\beta}$
(cf. \cite[Theorem 23.4]{Rockafellar:70}) and since $M$ is open
we have $M=\Int\left(\dom\hat{\beta}\right)$ and hence the convexity
of $\dom\hat{\beta}$ implies the convexity of $M$. Consequently,
(\ref{eq:Continuity}) immediately follows from \cite[Theorem 6.1]{Rockafellar:70}
and \cite[Corollary 7.5.1]{Rockafellar:70}. 

Clearly, by the definition of $\dom\hat{\beta}$, we have $\hat{\beta}\left(\alpha\right)=-\infty$
for $\alpha\not\in\overline{M}=\overline{\dom\hat{\beta}}$. The finiteness
of $\hat{\beta}$ on $\overline{M}$ follows from (\ref{eq:Continuity})
and from the fact that $0\leq\hat{\beta}\left(a\right)\leq\HD\left(\Lambda_{\Phi}\right)$
for all $a\in M$. This shows (\ref{eq:Mdom}).
\end{proof}
Now we are in the position to give a proof of the first part of Theorem
\ref{thm:MF}.

\subsection{Proof of Theorem \textmd{}\ref{thm:MF}}

\subsubsection{Proof of the first part of Theorem \ref{thm:MF}}

~

We split the proof of this theorem in two parts - upper bound and
lower bound.

For the upper bound we actually show a little more. Namely, for $\lambda\in\left(0,1\right)$
and $\alpha\in\overline{M}$ we consider $a_{\lambda}:=\alpha_{0}+\lambda\left(\alpha-\alpha_{0}\right)$,
where $\alpha_{0}:=\nabla\beta\left(0\right)$ is the unique maximum
of $\hat{\beta}$. For \[
\mathcal{G}_{a_{\lambda}}\left(v\right):=\left\{ \omega\in E_{A}^{\infty}:i\left(\omega_{1}\right)=v,\exists\ell\geq\lambda:\lim_{k\to\infty}\frac{S_{k}J\left(\omega\right)}{S_{k}I\left(\omega\right)}=\alpha_{0}+\ell\left(\alpha-\alpha_{0}\right)\right\} \]
we prove\begin{equation}
\HD\left(\mathcal{G}_{a_{\lambda}}\left(v\right)\right)\leq\hat{\beta}\left(a_{\lambda}\right).\label{eq:upperBound}\end{equation}
Because then we have by the monotonicity of the Hausdorff dimension
and (\ref{eq:Continuity})\[
\HD\left(\mathcal{F}_{\alpha}\left(v\right)\right)\leq\HD\left(\mathcal{G}_{a_{\lambda}}\left(v\right)\right)\leq\hat{\beta}\left(a_{\lambda}\right)\to\hat{\beta}\left(\alpha\right)\quad\textrm{for }\:\lambda\to1.\]
 To prove (\ref{eq:upperBound}) let us first define $b:\left(0,1\right)\to\R$,
$\lambda\mapsto\hat{\beta}\left(a_{\lambda}\right)$. Then $b'\left(\lambda\right)=-\left\langle t\left(a_{\lambda}\right),\alpha-\alpha_{0}\right\rangle $
is non-positive since $\alpha_{0}$ is the unique maximum of the strictly
concave function $\hat{\beta}$. Let $\mu_{t\left(a_{\lambda}\right)}$
denote the Gibbs measure for $\left\langle t\left(a_{\lambda}\right),J\right\rangle -\beta\left(t\left(a_{\lambda}\right)\right)I$.
Then by the above we have for $\epsilon>0$ and every $\omega\in\mathcal{G}_{a_{\lambda}}\left(v\right)$,
i.e. $\lim_{k\to\infty}\frac{S_{k}J\left(\omega\right)}{S_{k}I\left(\omega\right)}=\alpha_{0}+\ell\left(\alpha-\alpha_{0}\right)$
for some $\ell\geq\lambda$, that\begin{eqnarray*}
\mu_{t\left(a_{\lambda}\right)}\left(C_{n}\left(\omega\right)\right) & \gg & \exp\left(\left\langle t\left(a_{\lambda}\right),S_{n}J\left(\omega\right)\right\rangle -\beta\left(t\left(a_{\lambda}\right)\right)S_{n}I\left(\omega\right)\right)\\
 & = & \exp\left(-S_{n}I\left(\omega\right)\left(-\left\langle t\left(a_{\lambda}\right),\frac{S_{n}J\left(\omega\right)}{S_{n}I\left(\omega\right)}\right\rangle +\beta\left(t\left(a_{\lambda}\right)\right)\right)\right)\\
 & \gg & \exp\left(-S_{n}I\left(\omega\right)\left(\hat{\beta}\left(a_{\lambda}\right)+\epsilon-\left(\ell-\lambda\right)\left\langle t\left(a_{\lambda}\right),\alpha-\alpha_{0}\right\rangle \right)\right)\\
 & \gg & \left|\pi\left(C_{n}\left(\omega\right)\right)\right|^{\hat{\beta}\left(a_{\lambda}\right)+\epsilon}.\end{eqnarray*}
 Now consider a sequence of balls $\left(B\left(\pi\left(\omega\right),r_{n}\right)\right)_{n\in\N}$
with center in $\pi\left(\omega\right)\in X_{v}$ and radius $r_{n}:=\left|\pi\left(C_{n}\left(\omega\right)\right)\right|$,
where $\left|A\right|$ denotes the diameter of the set $A\subset\R^{D}$.
Then we have for all $\epsilon>0$ \begin{eqnarray*}
\mu_{t\left(a_{\lambda}\right)}\circ\pi^{-1}\left(B\left(\pi\left(\omega\right),r_{n}\right)\right) & \gg & \mu_{t\left(a_{\lambda}\right)}\left(C_{n}\left(\omega\right)\right)\gg r_{n}^{\hat{\beta}\left(a_{\lambda}\right)+\epsilon}.\end{eqnarray*}
Hence, by standard arguments from geometric measure theory (cf. \cite{Falconer:97,Mattila:95}),
the upper bound in (\ref{eq:upperBound}) follows.

For the lower bound we first consider $\alpha=\nabla\beta\left(t\right)\in M$.
Since $\int I\, d\mu_{t}<+\infty$ we have by \cite[Theorem 4.4.2]{MauldinUrbanski:03}
and the variational principle that \[
\HD\left(\mu_{t}\circ\pi^{-1}\right)=\frac{h_{\mu_{t}}}{\int I\, d\mu_{t}}=\frac{\beta\left(t\right)\int I\, d\mu_{t}-\left\langle t,\int J\, d\mu_{t}\right\rangle }{\int I\, d\mu_{t}}=\hat{\beta}\left(\alpha\right).\]
Since by Birkhoff's ergodic theorem we have $\mu_{t}\left(\pi^{-1}\left(\mathcal{F}_{\alpha}\right)\right)=1$
the above equality gives \[
\hat{\beta}\left(\alpha\right)=\HD\left(\mu_{t}\circ\pi^{-1}\right)\leq\HD\left(\mathcal{F}_{\alpha}\right).\]
The fact that by finite irreducibility $\HD\left(\mathcal{F}_{\alpha}\right)=\HD\left(\mathcal{F}_{\alpha}\left(v\right)\right)$
for all $v\in V$ finishes the proof of the lower bound for $\alpha\in M$.

For $\alpha\in\R^{d}\setminus\overline{M}$ we have on the one hand
that $\hat{\beta}\left(\alpha\right)=-\infty$ by (\ref{eq:Mdom})
and on the other hand by Theorem \ref{thm:Main2} that $\HD\left(\mathcal{F}_{\alpha}\left(v\right)\right)=\HD\left(\emptyset\right)=0$.
This proves the theorem for $\alpha\in\R^{d}\setminus\overline{M}$.

Before giving the proof of the remaining part of Theorem \ref{thm:MF}
we need the following proposition. Since from \cite{JenkinsonMauldinUrbanski:05}
we know that the entropy map is in our situation not upper semi-continuous
in general this proposition might be of some interest for itself. 

\begin{prop}
\label{pro:uscEntropy} We assume that $0\in M$ and let $\left(t_{j}\right)$
be a sequence in $\R^{d}$ with $\left|t_{j}\right|\to\infty$ such
that the sequence of Gibbs measures $\mu_{j}:=\mu_{t_{j}}$ converge
weakly to some $\mu\in\mathcal{M}\left(E_{A}^{\infty},\sigma\right)$.
Then $\mu$ is supported on a subshift of finite type over a finite
alphabet and we have\begin{equation}
\limsup_{j}h_{\mu_{j}}\leq h_{\mu}.\label{eq:uscEntropy}\end{equation}

\end{prop}
\begin{proof}
Since $0\in M$ guarantees that $\left|t_{j}\right|/\beta\left(t_{j}\right)$
stays bounded for $j\to\infty$ (cf. Subsection \ref{sub:Proof-of-L=3DoverbarM})
, $J$ is a bounded, and $I$ an unbounded function we find for $n\in\N$
and $D>0$ an $N\in\N$ such that for \[
\omega\in\left\{ x\in\Sigma_{A}^{n}:\exists m\in\left\{ 1,\ldots,n\right\} \, x_{m}\geq N\right\} =:\Sigma_{A}^{n}\left(N\right)\]
we have\begin{equation}
\left(\left\langle \frac{1}{\beta\left(t_{j}\right)}t_{j},S_{n}J\left(\omega\right)\right\rangle -S_{n}I\left(\omega\right)\right)<-D.\label{eq:lessThan0}\end{equation}
Using estimates from the proof of Theorem 2.3.3 in \cite{MauldinUrbanski:03}
(Gibbs property) one verifies that the constant $C$ in (\ref{eq:Gibbs})
is always less or equal to $\exp\left(K(\beta\left(t\right)+\left|t\right|)\right)$
for some positive constant $K$. Combining this fact, the Gibbs property
(\ref{eq:Gibbs}) and (\ref{eq:lessThan0}) then gives that \begin{equation}
\mu\left(C_{n}\left(\omega\right)\right)=0\;\textrm{for all }\:\omega\in\Sigma_{A}^{n}\left(N\right).\label{eq:FiniteSupport}\end{equation}
 Since $\mu$ is shift invariant it must be supported on a subshift
contained in the full shift $\left\{ 1,\ldots,N-1\right\} ^{\N}$. 

To show (\ref{eq:uscEntropy}) we set for $\nu\in\mathcal{M}\left(\Sigma_{A}^{\infty},\sigma\right)$
\[
H_{n}\left(\nu\right):=-\sum_{\omega\in\Sigma_{A}^{n}}\nu\left(C_{n}\left(\omega\right)\right)\log\nu\left(C_{n}\left(\omega\right)\right).\]
 It then suffices to verify that for all $n\in\N$ \begin{equation}
H_{n}\left(\mu_{j}\right)\to H_{n}\left(\mu\right)\quad\textrm{for }j\to\infty.\label{eq:approxConvergence}\end{equation}
 To see this note that $h_{\mu}=\lim\frac{1}{n}H_{n}\left(\mu\right)=\inf\frac{1}{n}H_{n}\left(\mu\right)$.
For $m,j\in\N$ we have\begin{eqnarray*}
h_{\mu_{j}} & = & \underbrace{\inf_{n}\frac{1}{n}H_{n}\left(\mu_{j}\right)-\frac{1}{m}H_{m}\left(\mu_{j}\right)}_{\leq0}+\frac{1}{m}H_{m}\left(\mu_{j}\right)-\frac{1}{m}H_{m}\left(\mu\right)+\frac{1}{m}H_{m}\left(\mu\right)\\
 & \leq & \underbrace{\frac{1}{m}H_{m}\left(\mu_{j}\right)-\frac{1}{m}H_{m}\left(\mu\right)}_{\to0,\;\textrm{as \,}j\to\infty}+\frac{1}{m}H_{m}\left(\mu\right)\end{eqnarray*}
implying $\limsup h_{\mu_{j}}\leq\frac{1}{m}H_{m}\left(\mu\right).$
Taking the infimum over $m\in\N$ gives (\ref{eq:uscEntropy}). Indeed
we have\begin{eqnarray*}
H_{n}\left(\mu_{j}\right) & = & -\sum_{\omega\in\Sigma_{A}^{n}}\mu_{j}\left(C_{n}\left(\omega\right)\right)\log\mu_{j}\left(C_{n}\left(\omega\right)\right)\\
 & = & -\sum_{\omega\in\Sigma_{A}^{n}\setminus\Sigma_{A}^{n}\left(N\right)}\mu_{j}\left(C_{n}\left(\omega\right)\right)\log\mu_{j}\left(C_{n}\left(\omega\right)\right)\\
 &  & \quad\quad\quad\,\qquad\;\qquad\;\qquad\;-\sum_{\omega\in\Sigma_{A}^{n}\left(N\right)}\mu_{j}\left(C_{n}\left(\omega\right)\right)\log\mu_{j}\left(C_{n}\left(\omega\right)\right).\end{eqnarray*}
The first sum has only finitely many summands and will therefore converge
to the sum $-\sum_{\omega\in\Sigma_{A}^{n}\setminus\Sigma_{A}^{n}\left(N\right)}\mu\left(C_{n}\left(\omega\right)\right)\log\mu\left(C_{n}\left(\omega\right)\right)$
as $j$ tends to infinity. Using the Gibbs property with constant
$C=\exp\left(K(\beta\left(t\right)+\left|t\right|)\right)$ one finds
that the second sum is summable and dominated by \begin{eqnarray*}
 &  & \!\!\!\!\!\!\exp\left(K\left(\left|t_{j}\right|+\beta\left(t_{j}\right)\right)\right)\sum_{\omega\in\Sigma_{A}^{n}\left(N\right)}\left(S_{n}\left(-\left\langle t_{j},J\right\rangle +\beta\left(t_{j}\right)I\right)\left(\omega\right)+K\left(\left|t_{j}\right|+\beta\left(t_{j}\right)\right)\right)\\
 &  & \qquad\qquad\qquad\qquad\qquad\qquad\qquad\qquad\;\;\;\qquad\times\exp\left(S_{n}\left(\left\langle t_{j},J\right\rangle -\beta\left(t_{j}\right)I\right)\left(\omega\right)\right)\\
\end{eqnarray*}
 The inequality in (\ref{eq:lessThan0}) guarantees that for $N$
sufficiently large this upper bound converges to $0$ as $j\to\infty$.
Hence, using (\ref{eq:FiniteSupport}) we conclude $H_{n}\left(\mu_{j}\right)\to H_{n}\left(\mu\right)$
for $j\to\infty$. 
\end{proof}

\subsubsection{Proof of the second part of Theorem \ref{thm:MF}}

~

Now we consider the special situation in which $0\in M$ and $\alpha\in\partial M$.
In Subsection \ref{sub:Proof-of-L=3DoverbarM} it has been shown that
there exists a sequence of Gibbs measures $\mu_{j}:=\mu_{t_{j}}$
converging weakly to some $\mu\in\mathcal{M}\left(E_{A}^{\infty},\sigma\right)$
such that $\left(\int I\, d\mu_{j}\right)^{-1}\int J\, d\mu_{j}=:\alpha_{j}$
lie on a line segment in $M$ for all $j\in\N$ and converges to $\alpha=\left(\int I\, d\mu\right)^{-1}\int J\, d\mu$.
Since $\mu$ is supported on a subshift of finite type over a finite
alphabet we may apply results from \cite[Appendix]{KesseboehmerStratmann:04a}
(which generalize results from \cite{Cajar}). Let $\mathcal{G}\left(\mu\right):=\left\{ \pi\left(x\right):x\in\Sigma_{A}^{\infty},\:\frac{1}{n}\sum_{j=0}^{n-1}\delta_{\sigma^{j}x}\stackrel{*}{\to}\mu\right\} $
be the set of the $\mu$ generic points. Then it has been shown in
\cite[Appendix]{KesseboehmerStratmann:04a} that there exists a Borel
probability measure $m$ on $E_{A}^{\infty}$ and a Borel set $M\subset E_{A}^{\infty}$
such that $\pi\left(M\right)\subset\mathcal{G}\left(\mu\right)$,
$m\left(M\right)=1$ and for all $x\in M$ we have \[
\liminf_{n\to\infty}\frac{\log\left(m\left(C_{n}\left(x\right)\right)\right)}{\log\left|\pi\left(C_{n}\left(x\right)\right)\right|}=\frac{h_{\mu}}{\int I\, d\mu}\:\textrm{and \,}\lim_{n}\frac{\log\left|\pi\left(C_{n}\left(x\right)\right)\right|}{-n}=\int I\, d\mu.\]
Now we argue similarly as in the proof of Theorem 4.4.2 in \cite{MauldinUrbanski:03}
to conclude that $\HD\left(M\right)\geq\left(\int I\, d\mu\right)^{-1}h_{\mu}$.
Fix $\epsilon>0$ small enough. Then by Egorov's Theorem, there exists
a Borel set $M'\subset M$ and $n_{0}\in\N$ such that $m\left(M'\right)>0$
and for all $n\geq n_{0}$ and $x\in M'$ we have\[
\frac{\log\left(m\left(C_{n}\left(x\right)\right)\right)}{\log\left|\pi\left(C_{n}\left(x\right)\right)\right|}\geq\frac{h_{\mu}}{\int I\, d\mu}-\epsilon\;\textrm{and }\:\frac{\log\left|\pi\left(C_{n}\left(x\right)\right)\right|}{-n}\leq\int I\, d\mu+\epsilon.\]
This gives for all $n\geq n_{0},$ $x\in M'$\[
m\left(C_{n}\left(x\right)\right)\leq\left|\pi\left(C_{n}\left(x\right)\right)\right|^{\frac{h_{\mu}}{\int I\, d\mu}-\epsilon}\:\textrm{and }\:\e^{-n\left(\int I\, d\mu+\epsilon\right)}\leq\left|\pi\left(C_{n}\left(x\right)\right)\right|\leq\e^{-n\left(\int I\, d\mu-\epsilon\right)}.\]
We fix $0<r<\exp\left(-n_{0}\left(\int I\, d\mu+\epsilon\right)\right)$
and for $x\in M'$ let $n\left(x,r\right)$ be the least number $n$
such that $\left|\pi\left(C_{n+1}\left(x\right)\right)\right|<r$.
By the above estimates we have that $n\left(x,r\right)+1>n_{0}$ and
hence $n\left(x,r\right)\geq n_{0}$ and $\left|\pi\left(C_{n\left(x,r\right)}\left(x\right)\right)\right|\geq r$.
Lemma 4.2.6 in \cite{MauldinUrbanski:03} guarantees that there exists
a universal constant $L\geq1$ such that for every $x\in M'$ and
$0<r<\exp\left(-n_{0}\left(\int I\, d\mu+\epsilon\right)\right)$
there exist points $x_{1},\ldots,x_{k}$ with $k\leq L$ such that
$\pi\left(M'\cap B\left(\pi\left(x\right),r\right)\right)\subset\bigcup_{\ell=1}^{k}\pi\left(C_{n\left(x_{\ell},r\right)}\left(x_{\ell}\right)\right).$
For $m'=m|_{M'}$ the restriction of $m$ to the set $M'$ we now
have\begin{eqnarray*}
m'\!\circ\!\pi^{-1}\left(B\left(\pi\left(x\right),r\right)\right) & \leq & \sum_{\ell=1}^{k}m\left(C_{n\left(x_{\ell},r\right)}\left(x_{\ell}\right)\right)\leq\sum_{\ell=1}^{k}\left|\pi\left(C_{n\left(x_{\ell},r\right)}\left(x_{\ell}\right)\right)\right|^{\frac{h_{\mu}}{\int I\, d\mu}-\epsilon}\\
 & \leq & \sum_{\ell=1}^{k}\e^{\left(-n\left(x_{\ell},r\right)\left(\int I\, d\mu-\epsilon\right)\left(\frac{h_{\mu}}{\int I\, d\mu}-\epsilon\right)\right)}\\
 & = & \sum_{\ell=1}^{k}\left(\e^{\left(-\left(n\left(x_{\ell},r\right)+1\right)\left(\int I\, d\mu+\epsilon\right)\right)}\right)^{\frac{n\left(x_{\ell},r\right)\left(\int I\, d\mu-\epsilon\right)}{\left(n\left(x_{\ell},r\right)+1\right)\left(\int I\, d\mu+\epsilon\right)}\left(\frac{h_{\mu}}{\int I\, d\mu}-\epsilon\right)}\\
 & \leq & \sum_{\ell=1}^{k}\left|\pi\left(C_{n\left(x_{\ell},r\right)+1}\left(x_{\ell}\right)\right)\right|^{\frac{n\left(x_{\ell},r\right)\left(\int I\, d\mu-\epsilon\right)}{\left(n\left(x_{\ell},r\right)+1\right)\left(\int I\, d\mu+\epsilon\right)}\left(\frac{h_{\mu}}{\int I\, d\mu}-\epsilon\right)}\\
 & \leq & \sum_{\ell=1}^{k}r^{\frac{n\left(x_{\ell},r\right)\left(\int I\, d\mu-\epsilon\right)}{\left(n\left(x_{\ell},r\right)+1\right)\left(\int I\, d\mu+\epsilon\right)}\left(\frac{h_{\mu}}{\int I\, d\mu}-\epsilon\right)}\\
 & \leq & L\cdot r^{\frac{\left(\int I\, d\mu-\epsilon\right)}{\left(\int I\, d\mu+\epsilon\right)}\left(\frac{h_{\mu}}{\int I\, d\mu}-2\epsilon\right)}\end{eqnarray*}
where the last inequality holds if we choose $n_{0}>\epsilon^{-1}\left(\frac{h_{\mu}}{\int I\, d\mu}-2\epsilon\right)$.
By the mass distribution principle this shows that $\frac{h_{\mu}}{\int I\, d\mu}\leq\HD\left(\pi\left(M'\right)\right)$.
Since $M'\subset M\subset\mathcal{G}\left(\mu\right)\subset\mathcal{F}_{\alpha}$
it follows that \[
\frac{h_{\mu}}{\int I\, d\mu}\leq\HD\left(\mathcal{F}_{\alpha}\right).\]
 By (\ref{eq:Continuity}) and Proposition \ref{eq:uscEntropy} we
therefore have\[
\hat{\beta}\left(\alpha\right)=\lim_{k\to\infty}\hat{\beta}\left(\alpha_{k}\right)=\lim_{k\to\infty}\frac{h_{\mu_{k}}}{\int I\, d\mu_{k}}\leq\frac{h_{\mu}}{\int I\, d\mu}\leq\HD\left(\mathcal{F}_{\alpha}\right).\]
This gives the lower bound for the Hausdorff dimension of $\mathcal{F}_{\alpha}$
also for $\alpha\in\partial M$.

\section*{Acknowledgment}

The first author was supported in part by Zentrale Forschungsförderung
Universität Bremen (Grant No. 03/106/2) and would also like to thank
the University of North Texas for its kind hospitality. The second
author was supported in part by the NSF Grant DMS 0400481.

\end{document}